\documentclass[11pt]{article}

\usepackage{amsbsy}
\usepackage{amssymb}
\usepackage{amsthm}
\usepackage{amsmath}
\usepackage{fullpage,hyperref,mathrsfs,txfonts,epsfig,graphicx,graphics}

\begin{document}

\title{The Euclidean distortion of the lamplighter group}
\author{Tim Austin\footnote{This work was conducted while T.
Austin was visiting the Courant Institute of Mathematical Sciences,
New York University.}\\UCLA\\{\tt timaustin@math.ucla.edu} \and
Assaf Naor\footnote{Research supported by NSF grants CCF-0635078 and
DMS-0528387.}\\Courant Institute\\{\tt naor@cims.nyu.edu} \and Alain
Valette\\Universit\'e de Neuch\^{a}tel\\{\tt
alain.valette@unine.ch}}
\date{}

\maketitle


\newenvironment{nmath}{\begin{center}\begin{math}}{\end{math}\end{center}}

\newtheorem{thm}{Theorem}[section]
\newtheorem{lem}[thm]{Lemma}
\newtheorem{prop}[thm]{Proposition}
\newtheorem{cor}[thm]{Corollary}
\newtheorem{conj}[thm]{Conjecture}
\newtheorem{dfn}[thm]{Definition}
\newtheorem{prob}[thm]{Problem}
\newtheorem{ques}[thm]{Question}


\newcommand{\A}{\mathcal{A}}
\newcommand{\B}{\mathcal{B}}
\newcommand{\K}{\mathcal{K}}
\renewcommand{\H}{\mathcal{H}}
\renewcommand{\Pr}{\mathrm{Pr}}
\newcommand{\s}{\sigma}
\renewcommand{\P}{\mathcal{P}}
\renewcommand{\O}{\Omega}
\renewcommand{\S}{\Sigma}
\newcommand{\T}{\mathrm{T}}
\newcommand{\co}{\mathrm{co}}
\newcommand{\e}{\mathrm{e}}
\newcommand{\im}{\mathrm{i}}
\renewcommand{\l}{\lambda}
\newcommand{\U}{\mathcal{U}}
\newcommand{\calH}{\mathcal{H}}
\newcommand{\G}{\Gamma}
\newcommand{\g}{\gamma}
\renewcommand{\L}{\Lambda}
\newcommand{\hcf}{\mathrm{hcf}}
\newcommand{\F}{\mathcal{F}}
\renewcommand{\a}{\alpha}
\newcommand{\bbN}{\mathbb{N}}
\newcommand{\bbR}{\mathbb{R}}
\newcommand{\bbZ}{\mathbb{Z}}
\newcommand{\bbC}{\mathbb{C}}
\newcommand{\C}{\mathbb{C}}
\newcommand{\bbE}{\mathbb{E}}
\newcommand{\bwr}{\boldsymbol{\wr}}

\newcommand{\bb}[1]{\mathbb{#1}}
\renewcommand{\rm}[1]{\mathrm{#1}}
\renewcommand{\cal}[1]{\mathcal{#1}}

\newcommand{\fin}{\nolinebreak\hspace{\stretch{1}}$\lhd$}

\begin{abstract}
\noindent We show that the cyclic lamplighter group $C_2 \bwr C_n$
embeds into Hilbert space with distortion ${\rm O}\left(\sqrt{\log
n}\right)$. This matches the lower bound proved by Lee, Naor and
Peres in~\cite{LeeNaoPer}, answering a question posed in that paper.
Thus the Euclidean distortion of $C_2 \bwr C_n$ is
$\Theta\left(\sqrt{\log n}\right)$. Our embedding is constructed
explicitly in terms of the irreducible representations of the group.
Since the optimal Euclidean embedding of a finite group can always
be chosen to be equivariant, as shown by Aharoni, Maurey and
Mityagin~\cite{AhaMauMit} and by Gromov (see~\cite{deCTesVal}), such
representation-theoretic considerations suggest a general tool for
obtaining upper and lower bounds on Euclidean embeddings of finite
groups.
\end{abstract}




\section{Introduction}

Given a bi-Lipschitz map $f:X \hookrightarrow Y$ from one metric
space $(X,\rho_X)$ into another $(Y,\rho_Y)$,  the
\textbf{distortion} of $f$ is defined to be the product of the
greatest expansion under $f$ and that under its inverse:
\[\rm{dist}(f) \coloneqq \sup_{\substack{x,z \in X\\ x\neq z}}\frac{\rho_Y(f(x),f(z))}{\rho_X(x,z)}\cdot
\sup_{\substack{x,z \in X\\ x\neq z}}
\frac{\rho_X(x,z)}{\rho_Y(f(x),f(z))}.\] We now define the overall
\textbf{distortion of $X$ into $Y$} to be the infimal distortion
over all bi-Lipschitz $f:X \hookrightarrow Y$ (and take this to be
$+\infty$ if no such maps exist), and write it $c_Y(X)$. There are
various contexts in which either a particular domain space or a
particular target space is of interest; for example, the distortions
of many different spaces into the Banach spaces $Y = L_p$ have been
studied extensively (see \cite{Mat02} for a partial survey of this
area).  In this case we write $c_p(X)$ in place of $c_{L_p}(X)$. In
this paper we will be concerned with the case $p = 2$, and will
refer to the distortion $c_2(X)$ as the \textbf{Euclidean
distortion} of $X$.  We will usually denote Hilbert space by
$\calH$, and will assume throughout that it is complex.

We will study the Euclidean distortion of a particular parameterized
family of groups: the cyclic lamplighter groups. These are defined
to be the wreath products of the order-two cyclic group $C_2 =
\{0,1\}$ by the cyclic groups $C_n = \{0,1,\ldots,n-1\}$. In
general, the \textbf{wreath product} $L \bwr H$ of some group $L$ by
some other group $H$ is the semidirect product $L^H \rtimes H$,
where $H$ acts on $L^H$  by left multiplication of the coordinates.
Concretely, $L\bwr H$ is the set $L^H\times H$ equipped with the
multiplication
\[\big((x_h)_{h \in H},g\big)\cdot\big((y_h)_{h \in H},k\big)
\coloneqq \big((x_h\cdot y_{gh})_{h \in H},gk\big).\]

Thus, our object of study will be $G \coloneqq C_2\bwr C_n$. Notice
that in this case the discrete cube $C_2^{C_n}$ appearing in the
definition of $G$ can be interpreted as the family $\cal{P}C_n$ of
subsets of $C_n$ by identifying $x = (x_k)_{k \in C_n}$ with $\{j
\in C_n:\ x_j = 1\}$, so that the group operation within this cube
is now the symmetric difference. Henceforth we will abuse notation
and treat a point $x \in C_2^{C_n}$ as a subset. This $G$ is a
finite solvable group, and can be generated by the two elements
$(\{0\},0)$ and $(\emptyset,1)$; these then give rise to a
left-invariant word metric $\rho$ on the group.

In \cite{LeeNaoPer} it was shown by the method of Markov convexity
that (with this metric understood) $c_2(G) \gtrsim \sqrt{\log n}$.
As noted in~\cite{LeeNaoPer}, an alternative proof of this lower
bound follows from exhibiting a constant distortion embedding of a
complete binary tree of depth $\Theta(n)$ into $G$
(see~\cite{LPP96}), and then applying Bourgain's lower bound for the
Euclidean distortion of trees~\cite{Bourgain86}. Somewhat
surprisingly, this embedded tree is an asymptotically worst-case
obstruction to embedding the entire lamplighter group $G$ into
Hilbert space. Our main result is that the above lower bound is
tight up to universal constants, answering a question posed
in~\cite{LeeNaoPer}:

\begin{thm}\label{thm:main}
For each $n$ there is a bi-Lipschitz map $f:G \hookrightarrow \calH$
for which
\[\rho\big((x,j),(y,k)\big) \lesssim \|f(x,j) - f(y,k)\| \lesssim \sqrt{\log n}\cdot \rho\big((x,j),(y,k)\big)\]
for all $(x,j),(y,k) \in G$.
\end{thm}

We will construct an embedding of $G$ of essentially least possible
distortion of a very special type: we will first specify an
\emph{action} $\beta$ of $G$ on a Hilbert space $\calH$ by unitary
operators (i.e. a unitary representation), and then obtain the
embedding into $\calH$ itself by carefully choosing a suitable point
$v \in \calH$ and then mapping $(x,j) \in G$ to the image of $v$
under $\beta(x,j)$.  Hilbert space embeddings of groups constructed
in this way are referred to as \textbf{equivariant}.

Note that if $G$ is locally compact and Abelian, then any map $f:G
\hookrightarrow \calH$ can be  analyzed via its vector-space valued
Fourier transform. The Euclidean embeddings of various Abelian
groups and some associated discrete spaces have been successfully
studied in this way: consider, for example, the analyses of flat
Riemannian tori and of quotients of the Hamming cube under group
actions in~\cite{KhoNao}. However, upon moving to non-Abelian groups
a general framework for either proving good lower bounds on their
Euclidean distortion or for isolating their low-distortion Euclidean
embeddings is yet to emerge; in addition to our use of an analysis
of irreducible representations to find such an embedding for the
group of interest here, we discuss in Section \ref{sec:coda} a
result, due to Aharoni, Maurey and Mityagin~\cite{AhaMauMit} in the
case of Abelian groups  and to Gromov (see~\cite{deCTesVal}) in the
case of general amenable groups, according to which equivariant
embeddings must always appear among those with minimal distortion.
We finish with some applications of this basic fact and  some open
problems.

\textbf{Remark on notation}\hspace{5pt} In addition to the Landau
notation ($\rm{o}$, $\rm{O}$, $\Omega$ and $\Theta$), in this paper
we will use $\approx$ and $\lesssim$, $\gtrsim$ to denote,
respectively, equality or the corresponding inequality up to some
universal positive multiplicative constant. We will also write
$\bbE[f(x)|x \in X]$ for the average of some function $f:X\to \bbC$
over a finite set $X$. \fin

\section{The embedding}

We will specify our embedding through an indexed family of
irreducible representations of the lamplighter group, together with
a vector in each of them. The direct sum of these representations
gives a single (fairly high-dimensional) representation of the
lamplighter group, together with the desired low-distortion
equivariant embedding into Hilbert space through the image of the
direct sum of these vectors.

\subsection{The lamplighter group and its
representations}\label{subs:warmup}


It seems helpful to recall the following heuristic description of
the lamplighter group with the aforementioned generators, if only
for the exposition of some of our later proofs. Consider a
collection of $n$ lamps indexed by $C_n$ (that is, say, positioned
equidistantly around a circular street), together with a lighter,
who walks along the street and either lights or douses lamps or
leaves them unchanged.

We now interpret a pair $(x,j)$ in $G$ as an operation on the whole
system of lamps and lighter: the lamps at those positions indexed by
the set $x \subseteq C_n$ will be changed (lighted if dark or
vice-versa), and the lighter will move to a position $j$ steps
further round the circle $C_n$. (Note that alternatively we could
think of $(x,j)$ as describing the state of the system with the
lamps at positions in $x$ illuminated and the lighter at position
$j$, but this intuition is a little less appropriate for
understanding the group law; of course, this `state' description of
$(x,j)$ simply arises by applying the `operation' $(x,j)$ to the
state with all lamps doused and the lighter initially at $0$.)

Given this description, we can think of the generator $(\{0\},0)$ as
the act of changing the lamp \emph{at the current location of the
lamplighter}, and similarly $(\emptyset,1)$ as the act of the
lamplighter moving one position to the next lamp around the circle.
Let us write $d_{C_n}$ for the obvious nearest-neighbour-graph
metric on the cyclic group $C_n$.

We shall use the following simple approximation for the word metric
$\rho$ on $G$.

\begin{lem}\label{lem:wordmetric}
The metric $\rho$ satisfies
\[\rho\big((x,j),(y,\ell)\big) \approx d_{C_n}(j,k) + \max_{k \in x\triangle y}\,(d_{C_n}(0,k)+1)\]
(where we interpret the maximum as $0$ if $x =y$).
\end{lem}

\textbf{Proof}\hspace{5pt} Since $\rho$ is an invariant metric it
suffices to show that for all $(x,j) \in G$
\[\rho\big((x,j),(\emptyset,0)\big) \approx d_{C_n}(0,j) + \max_{k\in
x}\,(d_{C_n}(0,k)+1).\] The $\rho$-distance of $(x,j)$ from
$(\emptyset,0)$ is the length of the shortest word in $(\{0\},0)$
and $(\emptyset,1)$ and their inverses that evaluates to $(x,j)$.
Certainly, such a word must contain at least $d_{C_n}(0,j)$ copies
of either $(\emptyset,1)$ or its inverse. Similarly, for any $k \in
x$, any word evaluating to $(x,j)$ must contain at least
$d_{C_n}(0,k)$ copies of the same generator, $(\emptyset,1)$, or its
inverse, since the lamplighter has to travel to position $k \in C_n$
in order to change the lamp at position $k$. In the latter case the
word must also contain at least one copy of $(\{0\},0)$ for the act
of changing that lamp. This proves that
$$
\rho\big((x,j),(\emptyset,0)\big) \ge d_{C_n}(0,j) + \max_{k \in
x}\,(d_{C_n}(0,k)+1).
$$

On the other hand, this reasoning shows at once that
$\rho\big((x,j),(x,0)\big)$ actually equals $d_{C_n}(0,j)$ for any
$x \in C_2^{C_n}$ (since no lamps need be lit or doused for this
journey of the lamplighter). In addition, for any $x \in C_2^{C_n}$,
the lamplighter can change all the lamps of $x$ by first traveling
to the furthest point of $x$ from $0$ on one side of $C_n$, lighting
the necessary lamps along the way, and then returning to the origin
and repeating this exercise on the other side.  This clearly takes
at most $6\max_{k \in x}\,(d_{C_n}(0,k) + 1)$ steps, and therefore
\begin{eqnarray*}
\rho\big((x,j),(\emptyset,0)\big) \leq\rho\big((x,j),(x,0)\big) +
\rho\big((x,0),(\emptyset,0)\big) \lesssim d_{C_n}(0,j) + \max_{k
\in x}\,(d_{C_n}(0,k)+1),
\end{eqnarray*}
as required. \qed

Next we recall some of the unitary representations of  $G$. Our list
consists of all the \emph{irreducible} representations when $n$ is
prime (these are found by the standard method of inducing
representations; see~\cite{FH91}). For composite $n$ some of these
representations break up further. However, we will only use members
of this list whole, and so will not trouble ourselves with the more
complicated decompositions for composite $n$. The representations of
interest fall naturally into two families:
\begin{itemize}
\item Some factor through the natural quotient mapping $C_2 \bwr C_n \twoheadrightarrow
C_n$ with kernel the normal subgroup $C_2^{C_n}\times\{0\}$, and
these are then given just by the (one-dimensional) irreducible
representations of $C_n$: for each $u \in \{0,1,\ldots,n-1\}$ we
obtain the character
\[\chi_u(x,j) \coloneqq \rm{e}^{2\pi\im u j/n}\rm{Id}_{\bbC},\]
where of course $\chi_0$ is just the trivial representation
$\pi_{\rm{triv}}$.
\item The remainder of our list corresponds to direct sums of
non-trivial one-dimensional irreducible representations of the cube
tied together by a permutation action of the lamplighter-motion
group $C_n$. Let $\alpha:C_n\to C_n$ denote the cyclic left shift
$\alpha(j)\coloneqq j-1$, and define for each $A \subseteq C_n$ the
{\bf Walsh function} $W_A:C_2^{C_n}\to \{-1,1\}$ by
$W_A(x)=(-1)^{|A\cap x|}$. For $A \notin \{\emptyset,C_n\}$  we
define the representation $\pi_A:C_2\bwr C_n \curvearrowright
\bbC^{C_n}$ by
\begin{eqnarray*}
\big(\pi_A(\emptyset,1)v\big)_k &\coloneqq& v_{k+1},\\
\big(\pi_A(x,0)v\big)_k &\coloneqq& W_A(\alpha^k(x))v_k =
(-1)^{|A\cap \alpha^k(x)|}v_k.
\end{eqnarray*}
For $A = C_n$ this is replaced by its more degenerate relative, the
one-dimensional representation \[\pi_{C_n}(x,j) \coloneqq
(-1)^{|x|}\rm{Id}_{\C}= W_{C_n}(\alpha^j(x))\rm{Id}_{\C}.\]
\end{itemize}

Note that there is a natural extension of the definition of $\pi_A$
to the case $A = \emptyset$:
\begin{eqnarray}\label{eq:emptyset}
\big(\pi_\emptyset(\emptyset,1)v\big)_k \coloneqq v_{k+1}\quad\quad
\hbox{and}\quad\quad \pi_\emptyset(x,0) \coloneqq
\rm{Id}_{\bbC^{C_n}};
\end{eqnarray}
this is given simply by composing the quotient $C_2\bwr C_n
\twoheadrightarrow C_n$ with the regular representation $C_n
\curvearrowright \bbC^{C_n}$, and as such it is isomorphic to the
direct sum of all the one-dimensional representations $\chi_u$ in
the first part of our list.

Before introducing our specific embedding, let us motivate the
construction by considering some generalities of the task of
constructing a low distortion equivariant embedding from these
ingredients. Suppose we have constructed an equivariant embedding
$f$ of $G$, expressed as
\[f(x,j) = \beta(x,j)v\]
for some unitary representation $\beta:G\curvearrowright \calH$ that
decomposes as
\[\beta(x,j) = \Big(\bigoplus_{u \in C_n} \bigoplus_{r=1}^{a_u}\chi_u(x,j)\Big)\oplus
\Big(\bigoplus_{\emptyset \neq A \subseteq
C_n}\bigoplus_{s=1}^{b_A}\pi_A(x,j)\Big),\] where $a_u, b_A
\in\mathbb N\cup\{0\}$ are multiplicities, and some vector
\[v = \Big(\bigoplus_{u \in C_n} \bigoplus_{r=1}^{a_u}v^{u,r}\Big)\oplus
\Big(\bigoplus_{\emptyset \neq A \subseteq
C_n}\bigoplus_{s=1}^{b_A}v^{A,s}\Big) \in \calH\] with $v^{u,r}$
(respectively $v^{A,s}$) lying in the $r^{\mathrm{th}}$
(respectively $s^{\rm{th}}$) subspace corresponding to a
subrepresentation $\chi_u$ (respectively $\pi_A$).

We can calculate a counterpart to Lemma~\ref{lem:wordmetric} for use
in our subsequent analysis:
\begin{eqnarray}\label{eq:parseval}
\|f(x,j) - f(\emptyset,0)\|^2 = \sum_{u \in
C_n}\sum_{r=1}^{a_u}\left|\rm{e}^{2\pi\rm{i}ju/n} -
1\right|^2\left|v^{u,r}\right|^2 + \sum_{\substack{A \subseteq
C_n\\A\neq \emptyset}}\sum_{s=1}^{b_A}\sum_{k \in
C_n}\left|W_A(\alpha^k(x))\cdot v^{A,s}_{k+j} - v^{A,s}_k\right|^2.
\end{eqnarray}

Suppose now that we know for $f$ the bounds
\[\rho\big((x,j),(y,k)\big) \leq \|f(x,j) - f(y,k)\| \leq D\rho\big((x,j),(y,k)\big)\]
(that is, $\rm{dist}(f) \le D$ and $f$ has been multiplied by a
scalar if necessary so that it is non-contractive).  Then one
natural approach to proving lower bounds for $D$ is to consider the
averages of the squared distances $\rho((x,j),(y,k))^2$ and
$\|f(x,j) - f(y,k)\|^2$ for $(x,j)$, $(y,k)$ in some subsets of $G$
for which the forms of the latter averages simplify in terms of our
orthogonal decomposition for $f$.  A lower bound for $D$ can be
obtained by comparing different such averages, as can some
information on which equivariant embeddings might be close to
attaining distortion $D$. Given the invariant nature of our
functions of interest, we can simplify our task further by setting
$(y,k) = (\emptyset,0)$ and averaging only over $(x,j)$.

For the sake of simplicity, we will consider this heuristic for $f$
in which the one-dimensional representations $\chi_u$ appear only
through the direct sum $\pi_\emptyset = \bigoplus_{u \in
C_n}\chi_u$, as described in \eqref{eq:emptyset}, and in which all
the multiplicities $b_A$ are $1$ (noting that if the vector $v^A$ is
$0$ then the representation $\pi_A$ effectively does not appear).

Suppose, then, that $H$ is some subset of $G$. Then
\begin{eqnarray}\label{eq:conditional}
\bbE\left[\rho\big((x,j),(\emptyset,0)\big)^2\,\big|\,(x,j) \in
H\right] &\leq& \bbE\left[\|f(x,j) -
f(\emptyset,0)\|^2\,\big|\,(x,j) \in H\right]\nonumber\\
&=&\nonumber \sum_{A \subseteq C_n}\bbE\left[\left\|\pi_A(x,j)v^A -
v^A\right\|^2\,\big|\,(x,j) \in H\right]\\&=& 2\sum_{A \subseteq
C_n}\left\|v^A\right\|^2 - 2\sum_{A \subseteq
C_n}\left\langle\bbE\left[\pi_A(x,j)\,\big|\,(x,j) \in
H\right]v^A,v^A\right\rangle.
\end{eqnarray}
This will be most helpful to us if we can arrange that for each
$A\subseteq C_n$ the expectation $\bbE\left[\pi_A(x,j)\,\big|\,(x,j)
\in H\right]$ takes a simple form on the copy of $\C^{C_n}$
corresponding to $\pi_A$.  This happens, for example, if $H$ is a
subcube of the canonical subgroup $C_2^{C_n}\times\{0\}$ of $G$. (In
fact, this can be fitted into a more general discussion of averages
over subgroups, but we postpone this to Subsection
\ref{subs:suffice}.) If $H = \{(x,0):\ x\subseteq B\}$ for some
$B\subseteq C_n$ then a straightforward calculation
reduces~\eqref{eq:conditional} to
$$
\bbE\left[\rho\big((x,0),(\emptyset,0)\big)^2\,\big|\,x\subseteq
B\right]\le 2\sum_{A\subseteq C_n}\sum_{\substack{k\in C_n\\ B\cap
\alpha^{-k}(A)\neq \emptyset}}\left|v^A_k\right|^2.
$$

On the other hand, we can apply the upper bound on $\|f(x,j) -
f(\emptyset,0)\|^2$ with $(x,j)$ one of the generators to see that
\[D^2 = D^2\rho\big((\{0\},0),(\emptyset,0)\big)^2 \geq \|f(\{0\},0) -
f(\emptyset,0)\|^2 = 2\sum_{\substack{A \subseteq C_n\\ 0\in
A}}\left\|v^{A}\right\|^2\] and
\begin{eqnarray*}
D^2= D^2\rho\big((\emptyset,1),(\emptyset,0)\big)^2 &\geq&
\|f(\emptyset,1) - f(\emptyset,0)\|^2= 2\sum_{A \subseteq
C_n}\sum_{k \in C_n}\left|v^{A}_{k+1} - v^{A}_k\right|^2.
\end{eqnarray*}

Finally, by actually estimating the expectation
$\bbE\left[\rho\big((x,0),(\emptyset,0)\big)^2\,\big|\,x\subseteq
B\right]$ using Lemma \ref{lem:wordmetric}, we can now use the above
two inequalities to give a lower bound for $D$ by comparing
\begin{eqnarray}\label{eq:our first sum}
\sum_{A\subseteq C_n}\sum_{\substack{k\in C_n\\ B\cap
\alpha^{-k}(A)\neq \emptyset}}\left|v^A_k\right|^2
\end{eqnarray}
against \begin{eqnarray}\label{eq:our second sum}\sum_{\substack{A \subseteq C_n\\
0\in A}}\left\|v^{A}\right\|^2\ \ \ \ \ \ \hbox{and}\ \ \ \ \ \
\sum_{A \subseteq C_n}\sum_{k \in C_n}\left|v^{A}_{k+1} -
v^{A}_k\right|^2\end{eqnarray} for different possible choices of
$v^A$.


Such a comparison might rely on the Poincar\'e inequality for the
discrete circle $C_n$, applied to the functions $v^{A}_\bullet$.
However, a careful examination now shows that playing with different
choices of $B\subseteq C_n$  does not give a non-trivial (that is to
say, growing in $n$) lower bound for $D$, even though we know from
the Markov convexity calculation that  $D \gtrsim \sqrt{\log n}$.

This very failure does, however, suggest that relatively
low-distortion embeddings might be found by looking for those $v^A$
that are close to saturating the Poincar\'e inequality for $C_n$.
For each $A\subseteq C_n$ this inequality bounds the overall average
squared difference
\[\frac{1}{n^2}\sum_{j,k \in C_n}\left|v^{A}_j - v^{A}_k\right|^2,\]
by a multiple of the local average
\[\frac{1}{n}\sum_{k \in C_n}\left|v^{A}_{k+1} - v^{A}_k\right|^2.\]
In general, the latter must be multiplied by $n^2$ to bound the
former, but this inequality is close to tight only if the function
$v^{A}_\bullet$ varies relatively slowly around the circle (that is,
if its Fourier transform is concentrated at low frequencies). One
finds that this near-saturation is necessary in order to obtain a
small distortion estimate from~\eqref{eq:our first sum}
and~\eqref{eq:our second sum} when $B = C_n$.  On the other hand,
for more general $B$ the resulting estimate can be kept small only
if we know that a positive proportion of the mass $\sum_{A\subseteq
C_n}\left\|v^A\right\|^2$ is contributed by sets that intersect $B$;
and this, in turn, requires that the distribution of the squared
norms $\|v^{A}\|^2$ be approximately invariant under rotations of
the sets $A$ and be spread roughly uniformly over sets $A$ of a
broad range of different sizes.

We have suppressed the calculations behind this discussion, as we
are presently trying only to be motivational.  One is led naturally
to consider sets $A \subseteq C_n$ that can be quite large, but are
not evenly distributed around $C_n$, so that there is some large arc
of $C_n$ away from $A$ into which we can concentrate most of the
$\ell^2(C_n)$-norm of a slowly-varying function $v^{A}_{\bullet}$.
In the next subsection we will construct an embedding from this
intuition, using all subsets $A$ that lie within some arc $I$ of the
circle $C_n$ of length $\lfloor n/3 \rfloor$.  Of course, we must
concede a distortion of least $\Omega(\sqrt{\log n})$ somewhere, and
it turns that this is manifested for the best possible choice of
$v^A$ in a slight shortfall from saturation of the Poincar\'e
inequality.

\subsection{The embedding itself}\label{subs:embedding}

The irreducible representations and corresponding vectors that we
will use will be indexed by pairs $(I,A)$ for $I$ an arc (i.e. a
connected subset) of $C_n$ of length $\lfloor n/3 \rfloor$ and $A
\subseteq I$. Let us write $\cal{I}$ for the family of such arcs, of
which there are $n$, and $\cal{P}I$ for the collection of subsets of
a given arc $I$. For each pair $(I,A)$, the corresponding
irreducible representation will simply be that indexed by $A$ in the
list of the previous subsection, retaining the convention that for
$A = \emptyset$ we identify $\pi_\emptyset$ with the
regular-quotient representation $C_2\bwr C_n \twoheadrightarrow C_n
\curvearrowright \bbC^{C_n},$ which is isomorphic to the direct sum
$\bigoplus_{u \in C_n}\chi_u$.

We  still need to specify the associated vector $v^{A,I}$.  We will
take this to depend only on $I$, defining $(v^{I}_{k})_{k \in C_n}$
by
\[v^{I}_{k} \coloneqq \left\{\begin{array}{ll}\eta&\ \ k \in I\\ \delta\sqrt{d_{C_n}(k,I)}&\ \ k\notin I.\end{array}\right.\]
This definition depends on the choice of the two parameters
$\delta$, $\eta$.  The analysis that follows below can be performed
by first allowing these to be free and then optimizing them; we
obtain
\begin{eqnarray*}
\eta \coloneqq\frac{1}{n2^{n/6}}\quad \mathrm{and} \quad \delta
\coloneqq \frac{1}{\sqrt{n}2^{n/6}}.
\end{eqnarray*}
  Another optimization is also implicit in our definition
of $v^{I}_{k}$: a priori, we could have taken $v^{I}_{k}$ to be of
the form $\delta \cdot d_{C_n}(k,I)^\alpha$ for $k \not\in I$ and
then optimized also over $\alpha > 0$.  This optimization does
naturally lead to the exponent $\alpha = \frac{1}{2}$: it turns out
that all other values of $\alpha$ give distortion following a
positive power law in $n$.

Note that this function $v^{I}_\bullet$ has the qualitative
properties suggested by our heuristic argument of the previous
subsection: it witnesses the small constant $1/n^2$ to within a
factor of $\log n$ for the Poincar\'e inequality on the circle
$C_n$, and has only a very small part of its $\ell^2(C_n)$-norm
inside the arc $I\supseteq A$. We have restricted ourselves to those
sets $A$ that can be contained in some arc $I$ of the circle
precisely so that in each summand with representation $\pi_A$ the
associated vector $v^{A,I}$ can be chosen to take small values on
$A$ but still be close to optimal for the Poincar\'e inequality.

Assembling the above, our overall embedding of $C_2 \bwr C_n$ is
given by
\begin{eqnarray*}
f(x,j) = \bigoplus_{I \in \cal{I}}\bigoplus_{A \in
\cal{P}I}\pi_A(x,j) v^I = \bigoplus_{I \in \cal{I}}\bigoplus_{A \in
\cal{P}I}\big(W_A(\alpha^k(x))\cdot v^{I}_{k+j}\big)_{k \in C_n}.
\end{eqnarray*}

We can now specialize the identity~\eqref{eq:parseval} to this data
(adjusting to our convention for $\pi_\emptyset$) and so compute:
\begin{multline}\label{eq:formula}
\|f(x,j) - f(\emptyset,0)\|^2 = \sum_{I \in \cal{I}}\sum_{A \in
\cal{P}I}\sum_{k \in C_n}\left|W_A(\alpha^k(x))v^{I}_{k+j} -
v^{I}_{k}\right|^2 \\\approx \sum_{I \in \cal{I}}\sum_{A \in
\cal{P}I}\Big(\sum_{k \in C_n}\left|v^{I}_{k+j} - v^{I}_{k}\right|^2
+ \sum_{k \in C_n}\boldsymbol{1}_{\{W_A(\alpha^k(x)) =
-1\}}\left|v^{I}_{k}\right|^2\Big),
\end{multline}
where in the second step we have used the additional fact that our
vectors $(v^{I}_{k})_{k \in C_n}$ have non-negative real entries, so
that
\begin{eqnarray*}
\left|W_A(\alpha^k(x))v^I_{k+j} - v^I_k\right|^2 =
\left|W_A(\alpha^k(x))(v^I_{k+j} - v^I_k) + (W_A(\alpha^k(x)) -
1)v^I_k\right|^2 \approx \left|v^I_{k+j} - v^I_k\right|^2 +
\boldsymbol{1}_{\{W_A(\alpha^k(x)) = -1\}}\left|v^{I}_{k}\right|^2.
\end{eqnarray*}

\textbf{Proof of Theorem \ref{thm:main}}\hspace{5pt} We prove the
upper and lower bounds on $\|f(x,j) - f(y,k)\|^2$ separately. Note
that since both this embedded distance and the original metric
$\rho$ are $G$-invariant it suffices to consider the case $(y,k) =
(\emptyset,0)$.

\textbf{Step 1: upper bound}\hspace{5pt} We wish to show that
\[\|f(x,j) - f(\emptyset,0)\| \lesssim \sqrt{\log n}\cdot \rho\big((x,j),(\emptyset,0)\big)\]
for all $(x,j) \in G$. Since $\rho$ is a word metric it suffices to
check this for $(x,j)$ equal to each of the two generators.

Suppose first that $(x,j) = (\{0\},0)$.  Then our
formula~\eqref{eq:formula} gives
\begin{eqnarray*}
\|f(\{0\},0) - f(\emptyset,0)\|^2\approx \sum_{I \in \cal{I}}\sum_{A
\in \cal{P}I}\sum_{k \in C_n}\boldsymbol{1}_{\{W_A(\alpha^k(\{0\}))
= -1\}}\left|v^{I}_k\right|^2= \sum_{I \in \cal{I}}\sum_{A \in
\cal{P}I}\sum_{k \in A}\left|v^{I}_{k}\right|^2= \eta^2n\sum_{A \in
\cal{P}I}|A| \approx \eta^2n^22^{n/3} = 1,
\end{eqnarray*}
owing to our choice of $\eta$.

Similarly, setting $(x,j) = (\emptyset,1)$, we obtain
\begin{eqnarray*}
\|f(\emptyset,1) - f(\emptyset,0)\|^2 \approx \sum_{I \in
\cal{I}}\sum_{A \in \cal{P}I}\Big(\sum_{k \in C_n}\left|v^I_{k+1} -
v^{I}_{k}\right|^2 + \sum_{k \in
C_n}\boldsymbol{1}_{\{W_A(\emptyset) =
-1\}}\left|v^{I}_{k}\right|^2\Big) = \sum_{I \in \cal{I}}\sum_{A \in
\cal{P}I}\sum_{k \in C_n}\left|v^{I}_{k+1} - v^{I}_k\right|^2.
\end{eqnarray*}
From our choice of $v^I$ we deduce that
\[\left|v^{I}_{k+1} - v^{I}_{k}\right| \approx \left\{\begin{array}{ll}0&\ \ \ \hbox{if}\ k,k+1 \in I\\
|\delta - \eta| \approx |\delta|&\ \ \ \hbox{if}\ |I \cap \{k,k+1\}|
= 1\\ \delta\frac{1}{\sqrt{d_{C_n}(k,I)}}&\ \ \ \hbox{if}\ k,k+1 \in
C_n\setminus I,\end{array}\right.\] and so the above sum can be
bounded by
\begin{eqnarray*}
\|f(\emptyset,1) - f(\emptyset,0)\|^2 &\lesssim& 2\delta^2\sum_{I
\in \cal{I}}\sum_{A \in \cal{P}I}\sum_{k = 1}^{\lfloor
n/3\rfloor}\left(\frac{1}{\sqrt{k}}\right)^2 =2\delta^2\sum_{I \in
\cal{I}}\sum_{A \in \cal{P}I}\sum_{k = 1}^{\lfloor
n/3\rfloor}\frac{1}{k} \approx 2\delta^2n2^{n/3}\log n \approx \log
n,
\end{eqnarray*}
owing to our choice of $\delta$.

Taking square roots and comparing these two estimates with the
approximation given by Lemma \ref{lem:wordmetric} completes the
check of both generators, and so also the proof of the upper bound;
note that these two checks already dictate our choice of $\eta$ and
$\delta$ up to $\rm{O}(\sqrt{\log n})$ and $\Omega(1/\sqrt{\log n})$
respectively.

\textbf{Step 2: lower bound}\hspace{5pt} We will obtain the lower
bound
\[\|f(x,j) - f(\emptyset,0)\| \gtrsim \rho\big((x,j),(\emptyset,0)\big)\]
by breaking the sum
\[\sum_{I \in \cal{I}}\sum_{A \in \cal{P}I}\Big(\sum_{k \in
C_n}\left|v^{I}_{k+j} - v^{I}_k\right|^2 + \sum_{k \in
C_n}\boldsymbol{1}_{\{W_A(\alpha^k(x)) =
-1\}}\left|v^{I}_{k}\right|^2\Big)\] into the two obvious subsums
and estimating these separately.

\textbf{Step 2.1: first sum}\hspace{5pt} We will use a rather crude
estimate obtained by considering various ranges of possible values
of $d_{C_n}(0,j)$ and for each of them summing over only a certain
range of $k$; this will be enough to obtain the lower bound we seek.

Observe from the definition of $v^{I}_k$ that if $d_{C_n}(0,j) \leq
d_{C_n}(k,I) \leq n/3 - d_{C_n}(0,j)$ then
\[\left|v^{I}_{k+j} - v^{I}_k\right| \gtrsim
\delta\frac{d_{C_n}(0,j)}{\sqrt{d_{C_n}(k,I)}}.\]

Suppose first that $d_{C_n}(0,j) \leq n/100$; then taking only those
$k$ in the above range  gives the lower bound
\begin{multline*}
\sum_{I \in \cal{I}}\sum_{A \in \cal{P}I}\sum_{k \in
C_n}\left|v^{I}_{k+j} - v^{I}_k\right|^2 \gtrsim \delta^2\sum_{I \in
\cal{I}}\sum_{A \in \cal{P}I}\left(\sum_{k:\,d_{C_n}(0,j) \leq
d_{C_n}(k,I) \leq n/3-d_{C_n}(0,j)}
\left(\frac{d_{C_n}(0,j)}{\sqrt{d_{C_n}(k,I)}}\right)^2\right)\\
\geq \delta^2d_{C_n}(0,j)^2\left(n2^{\lfloor n/3 \rfloor -
1}\right)\sum_{k = \lceil n/100 \rceil}^{\lfloor n/3 - n/100
\rfloor} \frac{1}{k} \approx
\left(\delta^2n2^{n/3}\right)d_{C_n}(0,j)^2 \approx d_{C_n}(0,j)^2,
\end{multline*}
recalling our choice of $\delta$.

On the other hand, if $d_{C_n}(0,j) > n/100$, then for those two
arcs $J_1$ and $J_2$ of points $k \in C_n$ satisfying $0 <
d_{C_n}(k,I) < n/1000$, at least one of them, say $J_1$, is such
that $d_{C_n}(k+j,I) \geq 10 d_{C_n}(k,I)$ for all $k \in J_1$.
These $k \in J_1$ therefore satisfy also
\[\left|v^{I}_{k+j} - v^{I}_k\right| \gtrsim \sqrt{n/1000}.\]
Therefore, taking instead the sum over $J_1$ in the above estimate,
we have
\begin{eqnarray*}
\sum_{I \in \cal{I}}\sum_{A \in \cal{P}I}\sum_{k \in
C_n}\left|v^{I}_{k+j} - v^{I}_k\right|^2 \gtrsim \delta^2\sum_{I \in
\cal{I}}\sum_{A \in \cal{P}I}\sum_{k\in
J_1}\left(\sqrt{n/1000}\right)^2\gtrsim
\left(\delta^2n2^{n/3}\right)(n/1000)^2 \gtrsim d_{C_n}(0,j)^2.
\end{eqnarray*}

In either case, we obtain
\[\sum_{I \in \cal{I}}\sum_{A \in \cal{P}I}\sum_{k \in
C_n}\left|v^{I}_{k+j} - v^{I}_k\right|^2\gtrsim d_{C_n}(0,j)^2.\]

\textbf{Step 2.2: second sum}\hspace{5pt} We now require a lower
bound on
\[\sum_{I \in \cal{I}}\sum_{A \in \cal{P}I}\sum_{k \in
C_n}\boldsymbol{1}_{\{|A \cap\alpha^k(x)|\
\rm{odd}\}}\left|v^{I}_k\right|^2= \sum_{I \in \cal{I}}\sum_{k \in
C_n}\left|v^{I}_k\right|^2\Big(\sum_{A \in
\cal{P}I}\boldsymbol{1}_{\{|A \cap\alpha^k(x)|\ \rm{odd}\}}\Big).\]

Note that for any non-empty subset $B$ of $C_n$ and for each $I \in
\cal{I}$, if we choose a subset $A$ of $I$ uniformly at random then
the probability that the intersection size $|A \cap B|$ is odd is
$1/2$ if $I \cap B \neq \emptyset$ and $0$ if $I \cap B =
\emptyset$. Indeed, choosing a subset $A$ of $I$ uniformly at random
and then considering $A\cap B$ simply generates a subset of $I \cap
B$ uniformly at random; but precisely half of these are odd unless
$I \cap B = \emptyset$, in which case they are all even. Applying
this reasoning with $B = \alpha^k(x)$, we conclude that
\[\sum_{A \in
\cal{P}I}\boldsymbol{1}_{\{|A \cap\alpha^k(x)|\ \rm{odd}\}} =
\frac{1}{2}|\cal{P}I|\boldsymbol{1}_{\{I \cap \alpha^k(x) \neq
\emptyset\}} = 2^{\lfloor n/3\rfloor - 1}\boldsymbol{1}_{\{I \cap
\alpha^k(x) \neq \emptyset\}},\] and so our sum of interest
simplifies to
\[2^{\lfloor
n/3\rfloor - 1}\sum_{I \in \cal{I}}\sum_{k \in
C_n}\boldsymbol{1}_{\{I \cap \alpha^k(x) \neq
\emptyset\}}\left|v^{I}_{k}\right|^2.\]

Suppose that $\ell \in x$ is a point of $x$ at a maximal distance
from $0$ in $C_n$. Then, in particular, $I \cap \alpha^k(x)
\supseteq I \cap \{\ell + k\}$ is nonempty for all $k \in
\alpha^{-\ell}(I)$, and so
\begin{eqnarray*}
2^{\lfloor n/3\rfloor - 1}\sum_{I \in \cal{I}}\sum_{k \in
C_n}\boldsymbol{1}_{\{I \cap \alpha^k(x) \neq
\emptyset\}}\left|v^{I}_{k}\right|^2 \geq 2^{\lfloor n/3\rfloor -
1}\sum_{I \in \cal{I}}\sum_{k \in
\alpha^{-\ell}(I)}\left|v^{I}_{k}\right|^2 .
\end{eqnarray*}
Therefore it will suffice to give a suitably strong lower bound for
$\sum_{k \in \alpha^{-\ell}(I)}\left|v^{I}_{k}\right|^2$. Moreover
we see from the rotational symmetry in our definition of $v^I$ that
this quantity is the same for all $I \in \cal{I}$. We may therefore
assume that in the natural labeling of $C_n$ as $\{1,2,\ldots,n\}$
the arc $I$ appears as an initial segment, and appealing to symmetry
further, we may replace $\ell$ by $-\ell$ and assume that $\ell \in
\{1,2,\ldots,\lfloor n/2\rfloor\}$. Given this, the terms appearing
in the desired sum are:
\begin{itemize}
\item terms equal to $\eta$ corresponding to $k \in I \cap
\alpha^{\ell}(I)$, and hence to $k \in \{\ell+1,\ldots,\lfloor
n/3\rfloor\}$;
\item the remaining terms $\delta,
\delta\sqrt{2},\delta\sqrt{3},\ldots,\delta\sqrt{\ell}$.
\end{itemize}
Squaring these and summing them therefore yields
\[\big(\lfloor n/3\rfloor - \ell\big)\eta^2 + \delta^2\sum_{t = 1}^{\ell}t \approx (n/3 - \ell)\eta^2 + \delta^2\ell^2,\]
and so overall
\[\sum_{I \in \cal{I}}\sum_{A \in \cal{P}I}\sum_{k \in
C_n}\boldsymbol{1}_{\{|A \cap\alpha^k(x)|\
\rm{odd}\}}\left|v^{I}_k\right|^2 \gtrsim n2^{n/3}\big((n/3 -
\ell)\eta^2 + \delta^2\ell^2\big) \gtrsim 1 + \ell^2,\] recalling
again our choices of $\delta$ and $\eta$.

\textbf{Completion of step 2}\hspace{5pt} Given the above estimates
for the first and second sum of our expression we deduce the lower
bound
\[\|f(x,j) - f(\emptyset,0)\|^2 \gtrsim d_{C_n}(0,j)^2 + 1 + \ell^2.\]
Recalling the choice of $\ell$, taking square roots and comparing
this with the expression of Lemma \ref{lem:wordmetric} completes the
proof. \qed


\section{Discussion and further questions}\label{sec:coda}

This section is composed of two parts.

In Subsection \ref{subs:suffice} we present the known result that
equivariant Euclidean embeddings of finite groups with invariant
metrics always appear among the embeddings of minimal distortion.
This justifies a reduction to the consideration of equivariant
embeddings of which we then give two applications.

In Subsection \ref{subs:questions} we discuss some further
questions.

\subsection{Equivariant embeddings suffice}\label{subs:suffice}

Unlike the special embedding of the lamplighter group constructed in
Subsection~\ref{subs:embedding}, a generic Hilbert space embedding
is certainly not equivariant. However, it turns out that searching
in this smaller class was, in a sense, guaranteed to work: for an
invariant metric on a finite group the restricted family of
equivariant embeddings must contain embeddings of distortion at
least as low as any other. This is the conclusion of
Lemma~\ref{lem:equivsuffice} below. For the same reason it suffices
to consider equivariant embeddings when proving Euclidean distortion
lower bounds for invariant metrics on finite groups. The formulation
 we give of Lemma~\ref{lem:equivsuffice} below is a simplified version for the case of finite
groups which we will use to investigate quantitative distortion
bounds. In the case of infinite Abelian groups this lemma was used
by Aharoni, Maurey and Mityagin~\cite{AhaMauMit} in their work on
uniform embeddings of Banach spaces into Hilbert space (see also
chapter 8 in the book~\cite{BenLin}). The lemma was discovered
independently by Gromov (unpublished) in the case of arbitrary
amenable groups, and was used by de Cornulier, Tessera and
Valette~\cite{deCTesVal} (in terms of Hilbert space valued cocycles)
to prove qualitative non-embeddability results for certain such
groups. Note that an analogous lemma holds for uniform embeddings
into Hilbert space, but for the sake of simplicity we present only
the bi-Lipschitz case.



\begin{lem}\label{lem:equivsuffice}
If a finite group $G$ with a left-invariant metric $\rho$ has a
Euclidean embedding $f$ such that
\begin{eqnarray}\label{eq:biLip}
\frac{1}{B}\|f(x) - f(y)\| \leq \rho(x,y) \leq A\|f(x) - f(y)\|
\end{eqnarray}
for all $x,y \in G$, then there is an equivariant embedding $g$ into
a Hilbert space $\calH$, say $g = \beta(\cdot)v$ for $v \in \calH$
and $\beta:G \curvearrowright\calH$, which satisfies the same
inequalities as in~\eqref{eq:biLip}.
\end{lem}

\textbf{Proof}\hspace{5pt} Define a positive semidefinite scalar
product on $\C^G$ by $K(\delta_x,\delta_y)\coloneqq \frac{1}{|G|}
\sum_{z\in G} \langle f(zx),f(zy)\rangle$. The required embedding
$g:G\to \C^G$ is simply given by $g(x)\coloneqq \delta_x$. Let
$\beta$ denote the left-regular representation of $G$ on $\C^G$.
Then $g(x)=\beta(x)\delta_e$, where $e$ is the identity element of
$G$. Now we compute that
\begin{eqnarray}\label{eq:summands}
\frac{\|g(x)-g(y)\|_K^2}{\rho(x,y)^2}=\frac{1}{|G|} \sum_{z\in G}
\left(\frac{\|f(zx)\|^2+\|f(zy)\|^2-2\langle
f(zx),f(zy)\rangle}{\rho(x,y)^2}\right)= \frac{1}{|G|} \sum_{z\in G}
\frac{\|f(zx)-f(zy)\|^2}{\rho(zx,zy)^2}.
\end{eqnarray}
By~\eqref{eq:biLip} each of the summands in~\eqref{eq:summands} lies
between $1/B^2$ and $A^2$, and hence so does the whole expression,
as required. It remains to note that $\|g(x)\|$ is independent of
$x\in G$, so that $\beta$ is a unitary representation with respect
to the scalar product $K$. \qed

Given this, we can now prove for arbitrary finite groups a
formalized version of the heuristic lower-bound analysis that was
presented in Subsection \ref{subs:warmup} to motivate the
construction of our embedding:

\begin{lem}\label{lem:average} Let $G$ be a finite group generated by $S\subseteq G$
and let $\rho$ be the corresponding word metric. Let
$\gamma_1\curvearrowright \H_1,\ldots, \gamma_t\curvearrowright
\H_t$ be the nontrivial irreducible representations of $G$. Then
there exist integers $a_1,\ldots, a_t\ge 0$ satisfying $\sum_{j=1}^t
a_j\dim(\H_j)\le |G|$ and sets of vectors
$\{v^{j,r}\}_{r=1}^{a_j}\subseteq \H_j$ for which
$$
c_2(G)\ge \sqrt{\frac{\sum_{x\in
G}\rho(x,e)^2}{2|G|}\cdot\frac{\sum_{s\in S}
\sum_{j=1}^t\sum_{r=1}^{a_j}\left\|\gamma_j(s)v^{j,r}-v^{j,r}\right\|^2}{
|S|\sum_{j=1}^t\sum_{r=1}^{a_j}\left\|v^{j,r}\right\|^2}}.
$$
\end{lem}

\textbf{Proof}\hspace{5pt} Assume that there exists a Euclidean
embedding $f$ satisfying~\eqref{eq:biLip}, and let $g$ be the
equivariant embedding from Lemma~\ref{lem:equivsuffice}. Note that
its dimension is at most $|G|$. We can write
$\beta=\bigoplus_{j=1}^t \gamma_j^{\oplus a_j}$, where $a_j\in
\mathbb N\cup\{0\}$ are multiplicities. Correspondingly we decompose
the vector $v$ from Lemma~\ref{lem:equivsuffice} as
$v=\bigoplus_{j=1}^t\bigoplus_{r=1}^{a_j} v^{j,r}$.  Then
\begin{eqnarray}\label{eq:compute}
\sum_{x\in G} \rho(x,e)^2&\le& \nonumber A^2 \sum_{x\in G}
\|g(x)-g(e)\|^2\\&=&\nonumber A^2\sum_{x\in G}
\sum_{j=1}^t\sum_{r=1}^{a_j}\left\|\gamma_j(x)v^{j,r}-v^{j,r}\right\|^2\\&=&
\nonumber A^2\sum_{x\in G} \sum_{j=1}^t\sum_{r=1}^{a_j}
\left(2\left\|v^{j,r}\right\|^2-2\left\langle\gamma_j(x)v^{j,r},v^{j,r}\right\rangle\right)\\
&=& \nonumber 2A^2|G|\sum_{j=1}^t\sum_{r=1}^{a_j}
\left\|v^{j,r}\right\|^2-2A^2\sum_{j=1}^t\sum_{r=1}^{a_j}\left\langle\Bigg(\sum_{x\in
G}\gamma_j(x)\Bigg)v^{j,r},v^{j,r}\right\rangle\\&=&2A^2|G|\sum_{j=1}^t\sum_{r=1}^{a_j}
\left\|v^{j,r}\right\|^2,
\end{eqnarray}
since $\sum_{x\in G}\gamma_j(x)=0$, by the irreducibility of
$\gamma_j$ (see~\cite{FH91}). On the other hand
$$
|S|=\sum_{s\in S}\rho(s,e)^2\ge \frac{1}{B^2}\sum_{s\in S}
\|g(s)-g(e)\|^2\ge \frac{1}{B^2}\sum_{s\in S}
\sum_{j=1}^t\sum_{r=1}^{a_j}\left\|\gamma_j(s)v^{j,r}-v^{j,r}\right\|^2.
$$
It follows that
$$
AB\ge \sqrt{\frac{\sum_{x\in
G}\rho(x,e)^2}{2|G|}\cdot\frac{\sum_{s\in S}
\sum_{j=1}^t\sum_{r=1}^{a_j}\left\|\gamma_j(s)v^{j,r}-v^{j,r}\right\|^2}{
|S|\sum_{j=1}^t\sum_{r=1}^{a_j}\left\|v^{j,r}\right\|^2}}.
$$
Infimizing over $AB$ yields the required result.
 \qed

\textbf{Remark}\hspace{5pt} We can obtain a larger family of lower
bounds for $c_2(G)$ by modifying the first part of the proof of
Lemma~\ref{lem:average} to the case of a sum over a subgroup $H$ of
$G$. However, this can lead to a more complicated expression owing
to the decomposition of the representations $\gamma_j$ into smaller
irreducible representations of $H$. Let $\pi_1\curvearrowright
\K_1,\ldots,\pi_m\curvearrowright \K_m$ be the irreducible
representations of $H$, where $\pi_1$ is the trivial representation
$\rm{Id}_{\K_1}$. Upon writing $\gamma_j$ as $\bigoplus_{\ell=1}^m
\pi_\ell^{\oplus b_{j,\ell}}$ and correspondingly $v^{j,r}$ as
$\bigoplus_{\ell=1}^m \bigoplus_{u=1}^{b_{j,\ell}}
v^{j,r}_{\ell,u}$, the sum $\sum_{x\in H}\gamma_j(x)$ equals
$\left(|H|\rm{Id}_{\K_1}^{\oplus b_{j,1}}\right)\oplus 0$. This
leads to a modification of~\eqref{eq:compute}, and thence to another
lower bound on the Euclidean distortion via the same argument.

In the case of the lamplighter group, this further decomposition
remains manageable, and it was by trying to approach equality in the
resulting lower bounds that we were led to the embedding of
Subsection~\ref{subs:embedding}.\fin

\textbf{Remark}\hspace{5pt} As indicated in Subsection
\ref{subs:warmup}, the lamplighter group $G = C_2 \bwr C_n$ has the
curious property that it does not embed into Hilbert space with
distortion bounded independent of $n$, but this nonembeddabability
is not detectable (in the sense of Lemma \ref{lem:average} and the
remark that follows it) by comparing the averages of the squared
group distances $\rho(x,y)^2$ and of the squared embedded distances
$\|f(x) - f(y)\|^2$ across subgroups of $G$ (the natural averages to
take) against the averages across local movements using the two
generators.  We find instead that for any given subgroup of $G$, $G$
itself has embeddings into Hilbert space that look good on average
across that subgroup, and `push' the bad distortion (which we know
must be at least $\Omega(\sqrt{\log n})$ somewhere) into some set of
pairs of point in the group that this average does not see. We
should stress that different subgroups may require slightly
different embeddings: the $\Omega(\sqrt{\log n})$ distortion of our
actual construction of Subsection \ref{subs:embedding}, for example,
can be detected by looking at averages across suitably-chosen
subgroups of $G$, while other embeddings, poorer overall, cannot be
detected by those subgroups. The point is that no small collection
of different subgroups reliably finds the distortion. This
conclusion follows from considering a number of variants of the
embedding of Subsection \ref{subs:embedding}; however, the necessary
calculations seem more lengthy than revealing and we will not
discuss them in detail here. Furthermore, one can also compute
easily given the methods of~\cite{NPSS06} that $G$, like Hilbert
space, does have Markov type $2$  (another averaging-based invariant
for metric spaces introduced by Ball in \cite{Bal}) with constant
independent of $n$, so that this also does not give an observable
obstruction to Hilbert space embeddings.

It results that both the minimal Euclidean distortion of $G$ and
also embeddings witnessing that distortion are hard to find using
standard averaging-based machinery.    Some quite delicate
averaging-based obstruction, such as the Markov convexity actually
used to study this group in \cite{LeeNaoPer}, is really necessary.
Furthermore, while that application of Markov convexity in
\cite{LeeNaoPer} does amount to the identification of a large
embedded tree in $G$, it is not at all clear a priori that the
minimal-distortion embeddings of this embedded tree already tell us
just how bad the Euclidean distortion of the whole group must be, or
how to attain that distortion.  It is somewhat surprising that this
invariant happens to give the correct growth rate of the Euclidean
distortion, and it might be interesting to ask whether Markov
convexity --- based, in this case, on the presence of large embedded
trees inside $G$ --- can be replaced by some averaging argument
using a different kind of substructure of $G$ to give the same lower
bound.\fin

We will finish our discussion of the consequences of Lemma
\ref{lem:equivsuffice} with a more concrete application.  Recall
that a metric space $(X,\rho)$ is of \textbf{negative type} if the
space $X$ with the square root metric $\sqrt{\rho}$ embeds
isometrically into Hilbert space.  The Goemans-Linial conjecture
asserted that any such metric also embeds with bounded distortion
into $L_1$.  This is now known to be false in general: see Khot and
Vishnoi~\cite{KV05}.  Their construction did not give an invariant
group metric; however, more recently Lee and Naor~\cite{LN06} have
shown that a particular invariant metric on the Heisenberg group is
also a counterexample to the Goemans-Linial conjecture, using a
result of Cheeger and Kleiner~\cite{CK06}.  On the other hand, the
following proposition shows that such counterexamples cannot arise
from Abelian groups subject to a restriction on the exponent of the
group (and we suspect that this restriction can be removed).

\begin{prop}\label{prop:gl} Let $(G,\rho)$ be a finite Abelian group
equipped with an invariant metric. Suppose that $2\le m\in \mathbb
N$ satisfies $mx=0$ for all $x\in G$. Let
$D=c_2\left(G,\sqrt{\rho}\right)$. Then $$c_1(G,\rho)\lesssim
D^4\log m$$ and for all $1<p<2$
$$c_p\left(G,\rho^{1/p}\right)\lesssim \frac{D^{4/p}}{p-1}.$$
\end{prop}

\textbf{Proof}\hspace{5pt} Let $\Gamma=\widehat G$ denote the dual
of $G$. By Lemma~\ref{lem:equivsuffice} $\left(G,\sqrt{\rho}\right)$
admits a distortion-$D$ embedding into Hilbert space which is
equivariant, and so which breaks into a direct sum of characters and
associated vectors. By rescaling it follows that there are
$\{a_\chi\}_{\chi\in \Gamma}\subseteq \mathbb R_+$ such that for all
$x\in G$,
\begin{eqnarray}\label{eq:dist}
\sum_{\chi\in \Gamma}a_\chi \left|1-\chi(x)\right|^2\le \rho(x,0)\le
D^2 \sum_{\chi\in \Gamma}a_\chi \left|1-\chi(x)\right|^2.
\end{eqnarray}
For every $x\in G$ and $k\ge 0$ denote
$$
A_k(x)\coloneqq \left\{\chi\in \Gamma:\ 2^{-k}<|\chi(x)-1|\le
2^{-k+1}\right\}.
$$
We also define $A_\infty(x)\coloneqq \{\chi\in \Gamma:\
\chi(x)=0\}$. Then  $\Gamma=A_\infty(x)\cup \bigcup_{k=0}^\infty
A_k(x)$ and this union is disjoint. Moreover,
$1=\chi(0)=\chi(mx)=\chi(x)^m$, so that $\chi(x)$ is an
$m^{\rm{th}}$ root of unity. Therefore if $\chi(x)\neq 0$ then
$|\chi(x)-1|\ge \left|{\rm e}^{2\pi \rm{i}/m}-1\right|\ge
\frac{1}{m}$. It follows that for finite $k>\log_2 m+1$ the set
$A_k(x)$ is empty.

It follows routinely from the definition of $A_k(x)$ that, firstly,
\begin{eqnarray}\label{eq:inA}
\sum_{\chi\in A_k(x)}a_\chi\left|1-\chi(x)\right|^2\ge
2^{-k}\sum_{\chi\in A_k(x)}a_\chi\left|1-\chi(x)\right|,
\end{eqnarray}
and, secondly, that for every $\chi\in A_k(x)$
\begin{eqnarray}\label{eq:twiddle}
\left|1-\chi(x)^{2^{k-1}}\right|\gtrsim 2^k|1-\chi(x)|.
\end{eqnarray}
Moreover, by the invariance of $\rho$ and the triangle inequality,
we know that for every $k\ge 1$,
\begin{eqnarray}\label{eq:triangle}
2^{k-1}\rho(x,0)\ge \rho\left(2^{k-1}x,0\right).
\end{eqnarray}
Therefore, for every $k\ge 1$ we have
\begin{multline}\label{eq:chain}
2^{k-1}\sum_{\chi\in \Gamma}
a_\chi|1-\chi(x)|^2\stackrel{\eqref{eq:dist}}{\ge}
\frac{2^{k-1}}{D^2}\rho(x,0)\stackrel{\eqref{eq:triangle}}{\ge}\frac{1}{D^2}\rho\left(2^{k-1}x,0\right)
\stackrel{\eqref{eq:dist}}{\ge} \frac{1}{D^2} \sum_{\chi\in \Gamma}
a_\chi\left|1-\chi(2^{k-1}x)\right|^2\\\ge
\frac{1}{D^2}\sum_{\chi\in
A_k(x)}a_\chi\left|1-\chi(x)^{2^{k-1}}\right|^2\stackrel{\eqref{eq:twiddle}}{\gtrsim}\frac{2^{2k}}{D^2}\sum_{\chi\in
A_k(x)}a_\chi\left|1-\chi(x)\right|^2\stackrel{\eqref{eq:inA}}{\gtrsim}
\frac{2^k}{D^2}\sum_{\chi\in A_k(x)}a_\chi\left|1-\chi(x)\right|.
\end{multline}
Thus (for $k\ge 1$ by~\eqref{eq:chain}, and trivially for $k=0$)
\begin{eqnarray}\label{eq:almost}
\sum_{\chi\in A_k(x)}a_\chi\left|1-\chi(x)\right|\lesssim
D^2\sum_{\chi\in \Gamma} a_\chi|1-\chi(x)|^2,
\end{eqnarray}
and so, combining the above,
\begin{multline}\label{eq:inL1}
\frac{\rho(x,0)}{D^2}\stackrel{\eqref{eq:dist}}{\le} \sum_{\chi\in
\Gamma} a_\chi |1-\chi(x)|^2\le 2\sum_{\chi\in \Gamma} a_\chi
|1-\chi(x)|=2\sum_{k\le \log_2 (2m)}\sum_{\chi\in A_k(x)} a_\chi
|1-\chi(x)|\\\stackrel{\eqref{eq:almost}}{\lesssim} D^2\log m
\sum_{\chi\in \Gamma}
a_\chi|1-\chi(x)|^2\stackrel{\eqref{eq:dist}}{\le} D^2\log m \cdot
\rho(x,0).
\end{multline}
Let $\mu$ be the measure on $\Gamma$ given by $\mu(\{\chi\})=a_\chi$
and consider the embedding $f:G\to L_1(\Gamma,\mu)$ defined by
$f(x)(\chi)\coloneqq\chi(x)$. Inequality~\eqref{eq:inL1} says
precisely that $\rm{dist}(f)\lesssim D^4\log m$. This completes the
proof of the first assertion of the proposition.

The proof of the second assertion is similar. Analogously
to~\eqref{eq:inA},
$$
\sum_{\chi\in A_k(x)}a_\chi\left|1-\chi(x)\right|^2\ge
2^{-k(2-p)}\sum_{\chi\in A_k(x)}a_\chi\left|1-\chi(x)\right|^p.
$$
Arguing as for~\eqref{eq:chain}
$$
\sum_{\chi\in A_k(x)}a_\chi\left|1-\chi(x)\right|^p\lesssim
D^22^{-(p-1)k}\sum_{\chi\in \Gamma} a_\chi|1-\chi(x)|^2.
$$
Hence
\begin{multline}\label{eq:inLp}
\frac{\rho(x,0)}{D^2}\le \sum_{\chi\in \Gamma} a_\chi
|1-\chi(x)|^2\le 2^{2-p}\sum_{\chi\in \Gamma} a_\chi
|1-\chi(x)|^p=2^{2-p}\sum_{k=0}^\infty\sum_{\chi\in A_k(x)} a_\chi
|1-\chi(x)|^p\\\lesssim D^2\sum_{k=0}^\infty 2^{-(p-1)k}
\sum_{\chi\in \Gamma} a_\chi|1-\chi(x)|^2\le
\frac{D^2}{1-2^{-(p-1)}} \cdot \rho(x,0)\lesssim \frac{D^2}{p-1}
\cdot \rho(x,0),
\end{multline}
and we may now take the same $f$ as above to give the required
embedding into $L_p(\Gamma,\mu)$.
 \qed

\textbf{Remark}\hspace{5pt} Proposition~\ref{prop:gl} implies, in
particular, that {\em any} invariant metric $\rho$ on the discrete
cube $\{0,1\}^d$ for which $c_2\left(\{0,1\}^d,\sqrt{\rho}\right)=D$
has also $c_1\left(\{0,1\}^d,\rho\right)\lesssim D^4$. It seems
likely that the fourth power is far from optimal. More
interestingly, for general finite Abelian groups we see no
compelling reason to believe that the factor of $\log m$ in
Proposition~\ref{prop:gl} is necessary. If it can be removed, this
would imply that no invariant metric on a finite Abelian group can
serve as a counterexample to the Goemans-Linial conjecture. (Note
that when $G$ is the cyclic group $C_m$ the factor $\log m$ can
already be improved to $\sqrt{\log m}\cdot \log \log m$ by the
general result of~\cite{ALN05}.) \fin


\subsection{Further questions}\label{subs:questions}

We speculate that the methods suggested by Lemma
\ref{lem:equivsuffice} can be used to give a fuller analysis of
minimal distortion Euclidean embeddings for various groups and
homogeneous spaces. For example, our experience with the lamplighter
group suggests that the methods of this paper may bear on other
wreath products $L \bwr H$; however, we should observe at once that
the behaviour of these can depend radically on the choice of
generators, even among those obtained by choosing generators for the
acting group $H$ and then including a suitable additional member of
$L^H$:

\begin{prop}\label{prop:wigderson} Let $S$ be a uniformly random
subset of $C_n$ of cardinality $|S|\ge 100\log n$, conditioned on
the event that it generates $C_n$ (which occurs asymptotically
almost surely). Consider the wreath product $C_2\bwr C_n$ equipped
with the word metric $\rho$ corresponding to the generating set
$\big(\{\emptyset\}\times S\big)\cup \{(\{0\},0)\}$. Then
$$
c_2\left(C_2\bwr C_n,\rho\right)\gtrsim \sqrt{n}
$$
asymptotically almost surely.
\end{prop}

\textbf{Proof}\hspace{5pt} Let $G$ be the Cayley graph of $C_n$ with
the generating set $S$. Then by~\cite{ALW01} the metric $\rho$ on
$C_2\bwr C_n$ is the shortest path metric on the zig-zag product of
the Hamming cube $C_2^{C_n}$ (with the standard graph structure) and
the graph $G$, which we denote by $H$ (we refer
to~\cite{ALW01,RVW02} for the definition and properties of the
zig-zag product). Let $\lambda_1$ be the normalized second
eigenvalue of $C_2^{C_n}$ and let $\lambda_2$ be the normalized
second eigenvalue of $G$. Then it is well known that
$\lambda_1=1-\frac{2}{n}$, and the Alon-Roichman
theorem~\cite{AloRoi} states that $\lambda_2$ is bounded away from
$1$ asymptotically almost surely. By Theorem 4.3 in~\cite{RVW02},
the normalized second eigenvalue of $H$, which we denote by
$\lambda$, is at most
$$
\frac12 \left(1-\lambda_2^2\right)\lambda_1+\frac12
\sqrt{\left(1-\lambda_2^2\right)^2\lambda_1^2+4\lambda_2^2}\le
1-\frac{\Omega(1)}{n},
$$
where the last estimate holds asymptotically almost surely. It
follows from a standard argument that for every $f:C_2 \bwr C_n\to
\cal H$,
$$
\frac{1}{|C_2 \bwr C_n|^2}\sum_{x,y\in C_2 \bwr C_n}
\|f(x)-f(y)\|^2\le \frac{2}{1-\lambda}\cdot
\frac{1}{|E(H)|}\sum_{xy\in E(H)} \|f(x)-f(y)\|^2.
$$
Hence, if $f$ satisfies $\rho(x,y)\le \|f(x)-f(y)\|\le D\rho(x,y)$
for every $x,y\in C_2 \bwr C_n$, then also
$$
\frac{1}{|C_2 \bwr C_n|^2}\sum_{x,y\in C_2 \bwr C_n}
\rho(x,y)^2\lesssim \frac{nD^2}{|E(H)|}\sum_{xy\in E(H)}
\rho(x,y)^2\le nD^2.
$$
Now observe that two randomly chosen points of the lamplighter group
differ in their first coordinate in a set of size $\Omega(n)$ with
probability bounded away from $0$. Therefore a positive proportion
of the terms in the left-hand side above are at least $\Omega(n^2)$,
since if two states $x$ and $y$ of the lamplighter group differ so
much, the lamplighter must change $\Omega(n)$ lamps to pass between
them, irrespective of the choice of generators in the ``movement
group" $C_n$. It follows that $D\gtrsim \sqrt{n}$, as required. \qed

It seems natural that the lower bound of $\Omega(\sqrt{n})$ obtained
above is tight, but we have not investigated this.

In spite of the possibility of a purely spectral argument as in
Proposition \ref{prop:wigderson}, we hope that a consideration of
equivariant embeddings and their decompositions into irreducible
representations may shed some light on other families of wreath
products, or other semidirect products.

Finally, we should note that our construction of Subsection
\ref{subs:embedding} clearly rests crucially on special properties
of Hilbert space embeddings, and so the following question remains
essentially untouched:

\begin{ques}
Does the infimal distortion $c_1(C_2\bwr C_n)$ of the lamplighter
group with the metric $\rho$ into the Banach space $L_1$ tend to
infinity with $n$?
\end{ques}




\bibliography{lamplighter}
\bibliographystyle{abbrv}

\parskip 0pt
\vspace{14pt}









\end{document}